\documentstyle[11pt,leqno]{article}
\input amssym.def
\input amssym.tex

\newtheorem{corr}{Corollary}[section]
\newtheorem{prop}{Proposition}[section]
\newtheorem{theo}{Theorem}[section]
\newtheorem{lemma}{Lemma}[section]
\newtheorem{claim}{Claim}[section]
\newtheorem{defi}{Definition}[section]
\newtheorem{rem}{Remark}[section]

\newtheorem{con}{Conjecture}[section]
\newcommand{\im}{{\rm im}}
\newcommand{\Alb}{{\rm Alb}}
\newcommand{\fgc}[2]{\hat{\pi}_{#1}(#2)}
\newcommand{\deck}{{\rm Deck}}

\begin{document}

\title{Nilpotent groups and  universal coverings of  smooth projective
varieties}
\author{L. Katzarkov \thanks{ The author was partially supported by A.P.
Sloan Dissertational Fellowship}}
\date{ }
\maketitle

\tableofcontents

\section{Introduction}

Characterizing the universal coverings of smooth projective
varieties is an old and hard question. Central to the subject is a conjecture
of Shafarevich according to which
the universal cover $\widetilde{X}$ of a smooth projective variety is
 holomorphically convex, meaning that for every  infinite sequence of points
without limit points in $\widetilde{X}$ there exists a holomorphic function
unbounded on this sequence.

\medskip

In this paper we try to study the universal covering of a smooth projective
variety $X$  whose fundamental group $\pi_{1}(X)$ admits an infinite
image homomorphism
\[ \rho : \pi_{1}(X) \longrightarrow L \]
into a complex linear algebraic group $L$. We will say that a nonramified
Galois covering $X' \rightarrow X$ corresponds to a representation
$\rho : \pi_{1}(X) \rightarrow L$ if its group of deck transformations is
$\im(\rho)$.

\begin{defi}
We call a representation $\rho : \pi_{1}(X) \rightarrow L$ linear,
reductive, solvable or nilpotent if the Zariski closure of its image is a
linear, reductive, solvable or nilpotent algebraic subgroup in $L$. We call
the corresponding covering linear, reductive, solvable or nilpotent
respectively.

The natural homomorphism $\pi_{1}(X,x) \rightarrow \fgc{{\rm uni}}{X,x}$ to
Malcev's pro-uni\-po\-tent completion will be called the Malcev representation
and the corresponding covering the Malcev covering.
\end{defi}

One may ask not only if the universal covering of $X$ is
holomorphically convex but also if some special intermediate coverings that
correspond to representations $\rho : \pi_{1}(X) \longrightarrow L$
are holomorphically convex.

\bigskip

In case $X$ is an algebraic surface and $\rho : \pi_{1}(X) \longrightarrow L$
is a reductive representation this question has been answered in \cite{KR}.
The author and M. Ramachandran proved there that if $X' \longrightarrow X$ is
a Galois covering of a smooth projective surface corresponding to a reductive
representation of $\pi_{1}(X)$ and such that  $\deck(X'/X)$ does not have two
ends, then  $X'$ is holomorphically convex. The proof is based on two major
developments in K\"{a}hler geometry that occured in the last decade. The first
is a correspondence, established through the work  of Hitchin \cite{HI},
Corlette \cite{C} and Simpson \cite{SC},
between Higgs bundles, representations of  the fundamental group  of a
smooth projective variety $\rho : \pi_{1}(X) \rightarrow G$
(here $G$ is a linear algebraic group over ${\Bbb C}$) and $\rho$  equivariant
harmonic maps from the universal covering of $X$ to the corresponding
symmetric space for $G$. This correspondence is called now - non-abelian
Hodge theory. The second is the theory of harmonic maps to buildings
developed by Gromov and Schoen
\cite{GS}.  This theory  gives the $p$-adic version of the theory of Higgs
bundles developed by Corlette, Hitchin and Simpson and can be thought as of
a $p$-adic non-abelian Hodge theory.

\medskip

These two ideas are used simultaneously in \cite{KR} in order to get more
information about $\pi_{1}(X)$. The proof  in \cite{KR} uses also some very
powerful ideas of Lasell, Ramachandran  \cite{BR} and Napier  \cite{N1}, which
can be
interpreted as a non-abelian strictness property. These ideas provide a bridge
and make the Nonabelian Hodge theory suitable for questions related to the
Shafarevich conjecture.

\medskip

In this paper we elaborate further on the  idea that the answer
to certain uniformization questions depends heavily  on the fundamental
group of the variety. We study the question if solvable or nilpotent coverings
$X' \rightarrow X$ are holomorphically convex for $X$ smooth projective
variety.

First we  prove  the following:

\begin{theo}
The Malcev covering of any  smooth projective $X$ is holomorphically
convex.
\end{theo}

As an immediate consequence of this statement we get:

\begin{theo} Let $X$ be a smooth projective variety with a virtually
nilpotent fundamental group. Then  the  Shafarevich conjecture is true for $X$.
\end{theo}
(Recall that a finitely generated group is nilpotent if its lower
central series has finitely many terms. A group is virtually nilpotent if it
has a finite index subgroup which is nilpotent.)

\bigskip

The proof of Theorem 1.1 uses the functorial  Mixed Hodge
Structure (MHS)  on $\pi_{1}(X)$ combined with  some new ideas of  J\'anos
Koll\'ar  from  \cite{K1} and  \cite{K2}.
At the end of section 2 we give a different proof  of Theorem 1.2, which
combined with the strictness property for the nonabelian Hodge theory  seems
to be a very promising idea ( see \cite{LM}).
Observe that what helps us prove Theorem 1.1 is the use of all nilpotent
representations of $\pi_{1}(X)$ at the same time.

We can ask even more basic question than the Shafarevich  conjecture:

\bigskip

\noindent
{\bf Question 1.} Are there any nonconstant holomorphic functions on the
universal covering $X$ of any smooth projective variety?

\bigskip

Clearly it is enough to restrict ourselves to the case when $\pi_{1}(X)$ is an
infinite group.

To study this question in bigger generality we add some more Hodge theoretic
tools - the  results of Arapura \cite{A}, Beauvile \cite{BE}, Green,
Lazarsfeld  \cite{GL} and Simpson \cite{SA}  about characterizing the
absolute sets in the moduli space of rank one local systems.  We also need
the following  variant of  the result of Arapura and Nori \cite{AN} saying
that linear solvable K\"{a}hler groups are nilpotent.

\begin{theo} Let $\Gamma$ be  a quotient of a  K\"{a}hler group $\pi_{1}(X)$
so that $\Gamma$ is a ${\Bbb Q}$-linear solvable group, then there are two
possibilities - either $\Gamma$ is
virtually nilpotent or $\pi_{1}(X)$ surjects onto the fundamental group of a
curve of genus bigger than zero.
\end{theo}

The above theorem gives a  way of constructing new examples of non-K\"{a}hler
groups.
In particular any group $\Gamma$ with infinite $H^{1}([\Gamma,\Gamma],
{\Bbb Q})$ possesing a solvable linear quotient defined over ${\Bbb
Q}$ that is not virtually nilpotent cannot be K\"{a}hler.

Unfortunately we could not prove a solvable variant of the theorem 1.1. The
maximum we were able to say is how much the solvable coverings differ from the
nilpotent ones. We show the following.

\begin{theo} If  $\Gamma$ is   a quotient of a  K\"{a}hler group $\pi_{1}(X)$
so that $\Gamma$ is a  complex linear solvable group, then there are two
possibilities - either $\Gamma$ is deformable to a
virtually nilpotent representation of $\pi_{1}(X)$  or $\pi_{1}(X)$ surjects
onto the fundamental group of a curve of genus bigger than zero.
\end{theo}

This  theorem gives a  way of constructing new examples of non-K\"{a}hler
groups.
In particular any group $\Gamma$ with infinite $H^{1}([\Gamma,\Gamma],
{\Bbb Q})$ possessing a solvable linear quotient defined over ${\Bbb
Q}$ that is not virtually nilpotent cannot be K\"{a}hler. In \cite{LM}  it is
proved that the linear covering are holomorphically convex for $X$ an algebraic
surface. Of course this implies a solvable variant of the theorem 1.1 for
algebraic surfaces. The above theorem implies immediately:

\begin{corr} Let $\rho : \pi_{1}(X) \longrightarrow S({\Bbb C})$ be a Zariski
dense  representation of the fundamental group of a smooth projective
variety $X$  to the complex points of
an affine solvable group defined over ${\Bbb Q}$. Then the image of
$\pi_{1}(X)$ in the Malcev completion of $\pi_{1}(X)$  is infinite.
\end{corr}

In particular this implies that the first Betti number of $X$ is nonzero so the
universal covering of $X$ $\widetilde{X}$ admits nonconstant holomorphic
functions. The above corollary can be proved of course in a different way too.

If we restrict ourselves to the case when $X$ is an algebraic surface we get a
stronger statement.

\begin{theo}  Let $ X $ be  a smooth projective surface  with an infinite
complex linear representation of its fundamental group. Then there exist
non-constant holomorphic functions on $\widetilde{X}$.
\end{theo}

In some sense the above theorem says that the  universal coverings are
different from arbitrary coverings. The well known example of Cousin  (see
e.g. \cite{N1}) gives a ${\Bbb Z}$-covering of the two  dimensional torus
which does not admit holomorphic functions.  The  Theorem 1.5 raises a natural
question:

\bigskip

\noindent
{\bf Question 2.} Are there examples of infinite  $\pi_{1}(X)$  without any
infinite linear representation?

\bigskip

There are known examples of groups with this properties, e.g. Higman's
four generator group. The question is if they can be fundamental groups.
Even more interesting question was asked by J. Koll\'ar and C. Simpson.

\bigskip

\noindent
{\bf Question 3.} Are there examples of infinite residually finite
$\pi_{1}(X)$  without any infinite linear representation?

\bigskip

As it was pointed out to me by S. Gersten the answer of this question is
positive if we are  looking for an arbitrary group not for $\pi_{1}(X)$.
There are the groups of Grigorchuk and Gupta-Sidki which are finitely
generated infinite torsion groups.  These groups
are known to be residually finite ( see e.g. \cite{BU}).

\medskip

A negative answer to this question could have a great impact on the answer to
Shafarevich conjecture  for residually finite groups (see \cite{LM},
\cite{LP}). From another side a recent paper of Bogomolov and the author
\cite{BL} shows that things can get quite exotic even for $\pi_{1}(X)$. In
some sense the examples constructed in \cite{BL} indicate that if the answer
of  {\bf Question 2.} is negative then the statement of Theorem 1.5
could be the best statement in such a  generality.

Theorem 1.5  and Corollary 1.1  suggest the following:

\begin{con}   Let $ X $ be  a smooth projective variety  with an infinite
linear representation of its fundamental group. Then there exist non-constant
holomorphic functions on $\widetilde{X}$.
\end{con}

All of this strongly suggests that Hodge theory has a lot to offer
in studying uniformization questions. We stop at the border line, before we
introduce the next level of Hodge theoretic considerations,  the theory of
Nonabelian Mixed Hodge Structures- a theory that is giving us  a way of
working with all linear representations at the same time to get maximal
information about $\pi_{1}(X)$. The first steps in this theory are done in
\cite{SL}, \cite{SIM1}, \cite{SIM2}, \cite{SIM3} and \cite{LP} and it is far
from being sufficiently developed.  In
any case it has fast  consequences even on a very primitive level.  Using
these very first steps we  prove in \cite{LM} the Shafarevich conjecture for
surfaces with linear fundamental groups. The same method  implies that the
coverings that correspond to any linear representation are holomorphically
convex. The proof  uses   basically only  the mixed Hodge structure on the
relative completion of $\pi_{1}(X)$ with respect to some complex variation of
mixed Hodge structures. Our feeling is that this is just the beginning.

\bigskip

\noindent
{\bf Acknowledgements:} I  would like to thank   A. Beilinson  F. Bogomolov, J.
Carlson,  K. Corlette,  R. Donagi, M. Gromov, S. Gersten, M. Larsen,  M. Nori,
T. Pantev, C. Simpson, D. Toledo and S. Weinberger for  the useful
conversations and H. Clemens,  P. Deligne,  R. Hain, J. Koll\'ar and M.
Ramachandran for teaching me all  ingredients of the technique used in this
paper. Special thanks to Professor J. Koll\'ar for inviting me to visit
University of Utah, where most of this work was done.

\section{The Malcev covering}

In this section we prove Theorem 1.1. and give some applications.

We start  with  some ideas of  J\'anos
Koll\'ar  from  \cite{K1} and  \cite{K2}.

In \cite{K1} Koll\'ar observed that the Shafarevich conjecture
is equivalent to:

1) There exists a normal variety ${\bf Sh}(X)$ and a proper
map with connected fibers   ${\bf Sh} : X
\longrightarrow {\bf Sh}(X) $, which contracts precisely  the   subvarieties
$Z$
 in $X$   with    the property that $ \im [\pi_{1}(Z')\longrightarrow
\pi_{1}(X)]$ is finite. Here $Z'$ denotes a desingularization of $Z$.

2) ${\bf Sh}(\widetilde{X} )$ is a Stein space. Here we denote by ${\bf
Sh}(\widetilde{X} )$ the Grauert- Remmert reduction of ${\bf Sh}(\widetilde{X}
)$. In our notations   ${\bf Sh}(X)={\bf Sh}(\widetilde{X} ) / \pi_{1}(X)$. The
action of $\pi_{1}(X)$ may have fixed points on ${\bf Sh}(\widetilde{X} )$ but
we can still take a quotient.

One can consider also a relative version of condition 1). Let $H
\triangleleft \pi_{1}(X)$ be a normal subgroup. We will say that
a subgroup $R \subset \pi_{1}(X)$ is almost contained in $H$ if
the intersection $R \cap H$ has finite index in $R$ and we will
write $R \lesssim H$.
We have the following condition.

\medskip

\begin{enumerate}
\item There exists a normal variety ${\bf Sh}^{H}(X)$ and a proper
map with connected fibers ${\bf Sh}^{H}: X \longrightarrow  {\bf Sh}^{H}(X),$
which contracts exactly the subvarieties  $Z$
in $X$ having the property that $\im[\pi_{1}(Z')\longrightarrow \pi_{1}(X)]
\lesssim H$. Again $Z'$ denotes a desingularization of $Z$.
The relative version of 2) is the following:
\item ${\bf Sh}^{H}(\widetilde{X} )$ is a Stein space. Here we denote by
${\bf Sh}^{H}(\widetilde{X} )$ the Grauert- Remmert reduction of
${\bf Sh}(\widetilde{X} )$. In our notations   ${\bf Sh}^{H}(X)
=  {\bf Sh}^{H}(\widetilde{X} ) / (\pi_{1}(X)/H)$.
\end{enumerate}
This was also independently observed by F. Campana in \cite{CM}.

\bigskip

Our approach is that if there is a natural candidate for ${\bf Sh}(X) $ it is
enough to check  condition 1) only for $Z$ - an algebraic curve. This certainly
is the case when  $\pi_{1}(X)$ is a nilpotent group. In the simplest case when
$\pi_{1}(X)$ is virtually abelian one uses for ${\bf Sh}(X) $ the Albanese
variety $\Alb (X)$.

It is clear  (see e,g, \cite{CT}) that for a smooth projective variety $X$
with   $\pi_{1}(X)$ an  infinite nilpotent group  the Albanese map:
\[ \Alb : X \longrightarrow \Alb(X) \]
has nontrivial image. In other words $\dim_{\Bbb{C}}(\im(\Alb))>0.$

Moreover if we denote by $S$ the Stein factorization of the Albanese map,
then this is a natural candidate for ${\bf Sh}(X) $ in case $\pi_{1}(X)$ is a
nilpotent group. Observe that the map

\[ X \longrightarrow S \]
contracts all subvarieties $Z$ with the property that $ \im [H_{1}(Z,
{{\Bbb Q}})\longrightarrow H_{1}(X, {{\Bbb Q}})]$ is  trivial.

Now using that $\pi_{1}(X)$ is a nilpotent group and the theory of Mixed Hodge
 Structures on its Malcev completion  we show that the fact that $ \im
[H_{1}(Z, {{\Bbb Q}})
\longrightarrow H_{1}(X, {{\Bbb Q}})]$ is  trivial is equivalent to the fact
that $ \im [\pi_{1}(Z)\longrightarrow \pi_{1}(X)]$ is finite for $Z$ an
algebraic curve. We finish the proof by reducing the argument for $Z$ of
arbitrary dimension to $Z$ an algebraic curve.

To prove Theorem 1.1 we need to show again that there is natural candidate for
 ${\bf Sh}^{H}(X) $,  where $H={\rm ker}(\rho : \pi_{1}(X) \longrightarrow
\fgc{{\rm uni}}{X,x}$ of the Malcev representation.
Again  this candidate is $S$ the Stein factorization of the Albanese map.
 At the end of section  we give a different proof  of Theorem 1.1, which is
basically spelling of the proof we have  given already in the language of
equivariant harmonic maps.

\subsection{Mixed Hodge Structure considerations}

In this subsection we explain why if  $\pi_{1}(X)$ is a nilpotent group  the
theory of Mixed Hodge Structures on it implies  that $\im [H_{1}(Z, {{\Bbb Q}})
\longrightarrow H_{1}(X, {{\Bbb Q}})]$ is  trivial is equivalent to the fact
that $ \im [\pi_{1}(Z)\longrightarrow \pi_{1}(X)]$ is finite for $Z$ an
algebraic curve. For some background one can look at \cite{D}, \cite{DG} or
\cite{H}.

For the proof of Theorem 1.1 we need to work with $X$ smooth but for
completeness in this section we will require only  the MHS on $H^{1}(X)$  is
of weights $> 0$.

\begin{lemma} If Z is a  compact nodal curve and
$f:Z \longrightarrow  X$ is a map to a  variety such that MHS on $H^{1}(X)$
is of weights $> 0$  then the map
\[ f_{* }: L(Z,x) \longrightarrow  L(X,f(x)) \]
is trivial if and only if the map
\[f^{*} : H^{1}(X,{\Bbb Q}) \longrightarrow H^{1}(Z,{\Bbb Q})\]
is trivial. Here $L(Z,x)$ and   $L(X,f(x))$ are the corresponding Lie algebras
 of the  unipotent completions    $\fgc{{\rm uni}}{Z,x}$  and $\fgc{{\rm uni}}{
X,f(x)}$ of
the fundamental groups $\pi_{1}(Z,x)$  and $\pi_{1}(X,f(x))$ respectively and
$x$ is a point in $Z$.
\end{lemma}
{\bf Proof.} Observe that  the map in unipotent completions determines and is
determined by a map on the corresponding Lie algebras:

\[L(Z,x) \longrightarrow  L(X,f(x)).\]
First let us  consider the case where $H_1(Z)$ is pure of weight $-1$. This
is the case when the dual graph of $Z$ is a tree. By a standard strictness
argument the weight filtration on  $L(Z,x)$ is its lower central
series, and the associated graded Lie algebra is generated by
$Gr_{-1} L(Z,x) = H_{1}(Z, {\Bbb Q})$.

Since
\[L(Z,x) \longrightarrow  L(X,f(x))\]
is a morphism of MHS, it is non-zero if and only if the
map
\[Gr L(Z,x) \longrightarrow Gr L(X,f(x))\]
on weight graded quotients is. Since
\[Gr_{-1} L(X,f(x)) = H_{1}(X,{\Bbb Q})/W_{-2},\]
and since $H_{1}(Z,{\Bbb Q}) \longrightarrow H_{1}(X,{\Bbb Q})$ is trivial, it
follows that $L(Z,x) \longrightarrow  L(X,f(x))$  is trivial.

\medskip

To prove the general case, we take a partial normalization
\[Z' \longrightarrow Z\]
with the property that $Z'$  is connected and such that $H^{1}(Z')$  is a pure
MHS  of weight 1.

This can be done as follows.
Take a maximal tree $T$  in the dual graph
of $Z$  and normalize only those double points corresponding to
edges not in $ T$. Then $H_{1}(Z) $ is pure MHS  of weight -1. The previous
argument implies  that
\[L(Z',x) \longrightarrow  L (X,f(x))\]
is trivial.

To complete the proof, note that we have an exact sequence
\[1 \longrightarrow  N \longrightarrow  \pi_{1}(Z,x) \longrightarrow
\pi_{1}(\Gamma,*) \longrightarrow  1,\]
where $\Gamma$ denotes the dual graph of $Z $ and $N$  is the normal subgroup
of $\pi_{1}(Z)$  generated by $\pi_{1}(Z',x)$. After passing to unipotent
completions, we obtain an exact sequence
\[0  \longrightarrow (L(Z'))\longrightarrow L(Z,x) \longrightarrow
L(\Gamma,*) \longrightarrow  0.\]
This is an exact sequence in the category of Malcev Lie
algebras with MHS. The ideal $(L(Z'))$ generated by $L(Z')$
is exactly $W_{-1} L(Z)$, so  the MHS induced
on $L(\Gamma,*)$  is pure of weight 0.

It follows that the homomorphism  $L(Z,x) \longrightarrow  L(X,f(x))$
induces a homomorphism
\[L(\Gamma,*) \longrightarrow L(X,f(x)).\]
This is a morphism of MHS of (0,0) type. It is injective
if and only if the map
\[L(\Gamma,*) = Gr L(\Gamma,*) \longrightarrow Gr L(X,f(x))\]
is also injective. Since $H_{1}(X)$ has weights $< 0$  and $L(\Gamma,*)$
has weight zero, it follows that
\[L(\Gamma,*) = Gr L(\Gamma,*) \longrightarrow Gr L(X,f(x))\]
 is zero.  This proves the statement in general.
Namely, we have that for any nodal curve (singular, reducible)  the map
\[ f_{* }: L(Z,x) \longrightarrow  L(X,f(x)) \]
is trivial if and only if the map
\[f^{*} : H^{1}(X,{\Bbb Q}) \longrightarrow H^{1}(Z,{\Bbb Q})\]
is trivial. \hfill $\Box$

\begin{lemma}  Let $X$ be a smooth projective variety with a nilpotent
fundamental group  $\pi_{1}(X)$. Then  for any algebraic curve $Z \subset X $
   the fact  $ \im [H_{1}(Z, {{\Bbb Q}})\longrightarrow
 H_{1}(X, {{\Bbb Q}})]$ is  trivial is equivalent to the fact that $ \im
[\pi_{1}(Z)\longrightarrow \pi_{1}(X)]$ is finite.

\end{lemma}
{\bf Proof.}  Since we
can always find a partial normalization $\widetilde{Z} \rightarrow Z$ with
$\widetilde{Z}$-nodal and $\pi_{1}(\widetilde{Z},\tilde{x}) \rightarrow
\pi_{1}(Z,x)$ surjective it follows from the previous lemma that the map
\[ f_{* }: L(Z,x) \longrightarrow  L(X,f(x)) \]
is the zero map.
Furthermore, if $\pi_{1}(X)$ is a torsion free nilpotent group then by
definition it embeds in  $\pi_{un}(X, f(x))$. It is easy to see that torsion
elements  of a nilpotent group generate a finite group and hence
\[\pi_{1}(X,f(x)) \longrightarrow \fgc{{\rm uni}}{X,f(x)}\]
is an embedding up to  torsion which proves the lemma.
\hfill
$\Box$

We have actually proved more:

\begin{lemma}  Let $X$ be a smooth projective variety and $\rho : \pi_{1}(X)
 \longrightarrow L(X,f(x))$ be the Malcev representation of  $\pi_{1}(X)$.
Then  for any algebraic curve $Z \subset X $    the fact  $ \im [H_{1}(Z,
{{\Bbb Q}})\longrightarrow
 H_{1}(X, {{\Bbb Q}})]$ is  trivial is equivalent to the fact that $ \im
[\pi_{1}(Z) \longrightarrow \pi_{1}(X) / H ]$ is finite. Here $H$ is the
kernel  of the Malcev representation.

\end{lemma}

\subsection{A reduction to the case of an algebraic curve}

In this section we show how to reduce the argument for $Z$ of arbitrary
dimension to $Z$ an algebraic curve.

\begin{lemma}  Let $F$ be a connected subvariety in $X$ then we can find a
curve $Z
\subset F$ such that $\pi_{1}(Z) $ surjects on $\pi_{1}(F) $.
\end{lemma}
{\bf Proof.} If $F$ is smooth variety the above lemma is just the Lefschetz
hyperplane section theorem. Let  $F= F_{1}+ \ldots + F_{i}$ be singular and
with many components of different dimension. Denote by $n$   the normalization
 $n:F' \longrightarrow F $  of $F$. In every  component of $F'$  after
additional desingularization we can find finitely many points $x_{k}, y_{k}$
such that  $n(x_{k})=n(y_{k})$ and  $\pi_{1}(F'/ x_{k}=y_{k})$ surjects onto
  $\pi_{1}(F)$.   The way to do that is to take the Whitney stratification of
$F$ and put the points $x_{k},y_{k}$ in every stratum in a way that all loops
that come from singularities pass through these points. Now following
\cite{GM}(ii, 1.1) we take    hypesurfaces with big degrees that pass through
 the points  $x_{k}, y_{k}$ and intersect  every component  of $F'$,  $F'_{l}$
 in a curve $Z_{l}$  such  that   $Z'= \cup Z_{l}$ and  $\pi_{1}(Z')$ surjects
 on $\pi_{1}(F') $. We make   $Z=n(Z')$. Observe that $Z$ might be singular
and have many components but it will be connected.
\hfill $\Box$

Now we are ready to finish the proof of Theorem 1.1. We start with the
Stein factorization of the Albanese map for $X$
\[ \Alb : X \longrightarrow S \longrightarrow \im(\Alb) \subset \Alb(X).\]
Denote by $S'$ the  fiber product of the universal covering  $\widetilde{\Alb(
X)}$ of $\Alb(X)$ and $S$ over $\Alb(X)$. By definition   the map
\[ S' \longrightarrow \widetilde{\Alb(X)} \]
is a covering  map  and since  $\widetilde{\Alb(X)}$ is a Stein manifold $S'$
is a Stein manifold as well. It follows from the  definition of the Albanese
morphism  that the fibers of the map
\[ \Alb : X \longrightarrow S \]
are all  subvarieties $F$ in $X$ for which the map
$H_{1}(F, {{\Bbb Q}})\longrightarrow H_{1}(X, {{\Bbb Q}})$ is  trivial. We
willshow that if  the fact that $H_{1}(F, {{\Bbb Q}})
\longrightarrow H_{1}(X, {{\Bbb Q}})$ is trivial implies that $ \im
[\pi_{1}(F)\longrightarrow \pi_{1}(X) / H ]$ is finite.

We have shown this in Lemma 2.3  when $F$ is an algebraic curve.
Now if $dim_{\Bbb{C}}(F)>1$ we apply Lemma 2.4 to find a curve $Z$ in $F$
such that $\pi_{1}(Z) $ surjects on $\pi_{1}(F) $. The argument of
Lemma 2.1.2  implies that  $\pi_{1}(F) $ goes to a finite group in
$\pi_{1}(X)/H $ since
$\pi_{1}(Z) $ goes to  afinite group  in  $\pi_{1}(X)/H $. Observe that the
curve $Z$ is also contained in the fiber $F$ of the map
\[ \Alb : X \longrightarrow S .\]
Therefore the map
$H_{1}(Z, {{\Bbb Q}})\longrightarrow H_{1}(X, {{\Bbb Q}})$ is  also trivial.
To finish the proof of Theorem 1.1 we need to observe that  $S$ satisfies
the conditions for being the Shafarevich variety of $X$, $S= {\bf Sh}^{H}(X)$.
Namely

1) There exists  a holomorphic map with connected fibers  $ X
\longrightarrow S $, which contracts only the   subvarieties  $Z$
 in $X$   with    the property that $ \im [\pi_{1}(Z)\longrightarrow
\pi_{1}(X)/H]$ is finite.

2)   $  {\bf Sh}^{H}(\widetilde{X})=S'$  is  a Stein space.

\hfill $\Box$

To prove    Theorem 1.2   we use  the same argument as above but $H$ is a
finite
group.

Actually we have shown more:

\begin{corr} Let $X$ be a smooth projective variety with a virtually
residually nilpotent
fundamental group. Then  the  Shafarevich conjecture is true for $X$.
\end{corr}

\subsection{Some examples}

In this subsections we give some examples and geometric applications of our
method. We start with the following  result that was also proved by Campana in
\cite{CM1}.

\begin{corr} Let  $X$ be  a smooth projective surface and  $\Gamma$  is the
image of $\pi_{1}(X)$  in $L(X,f(x))$. Let as before $S$  be  the Stein
factorization of the map $X \longrightarrow \im(\Alb(X))$. After taking an
etale finite covering $X'' \longrightarrow X$ the homomorphism $\pi_{1}(X'')
\longrightarrow \Gamma$ factors through the map  $\pi_{1}(S) \longrightarrow
\Gamma$.
\end{corr}
{\bf Proof.} According to \cite{K2} 4.8  after taking  some  etale finite
covering $X'' \longrightarrow X$, $\pi_{1}(X'')$  is the same as the
fundamental group of  $\pi_{1}(S)$. This follows from the fact that residually
nilpotent  groups  are  residually finite.

\hfill $\Box$

Nilpotent  K\"{a}hler groups were constructed by Sommese and Van de Ven
\cite{SV}, and  Campana \cite{CM1}. The construction goes as follows:

Start with a finite morphism from an abelian variety $A$  to ${\Bbb P}^{n}$.
Now take the preimage $X$ in $A$  of any abelian d-fold in  ${\Bbb P}^{n}$. A
double cover of $X$  has as fundamental group a nonsplit central extension of
an abelian group by ${\Bbb Z}$.

Let us following \cite{SV} give  more  explicit example. We start with a
four dimensional abelian variety $A$ and a finite  morphism $f$ to
${\Bbb P}^{4}$. Take the Mumford-Horrocks abelian surface $Z$  in ${\Bbb
P}^{4}$
 and pull it back to $A$. Let us call the new surface
$f^{-1}(Z)$. The following exact sequence was established in  \cite{SV}
\[\pi_{2}(A)\oplus \pi_{2}(Z) \longrightarrow  \pi_{2}({\Bbb P}^{4})
\longrightarrow \pi_{1}(f^{-1}(Z)) \longrightarrow \pi_{1}(A)\oplus
\pi_{1}(Z) \longrightarrow 0.\]
In our case this sequence reads as:
\[0 \longrightarrow {\Bbb Z} \longrightarrow \pi_{1}(f^{-1}(Z)) \longrightarrow
 {\Bbb Z}^{12} \longrightarrow 0\]
and  shows that $f^{-1}(Z)$ has a two steps nilpotent fundamental group.

Actually we know more. By   theorem of Arapura
and Nori \cite{AN}  all K\"{a}hler   linear solvable groups are
virtually nilpotent. So  we have the following:

\begin{corr} Let $X$ be a smooth projective variety  with a   linear solvable
 fundamental group. Then  the  universal  covering $\widetilde{X}$ is
holomorphically convex.
 \end{corr}

Now we will use the technique from Lemmas 2.3 and 2.4 to show that the
theorem of Arapura
and Nori \cite{AN} is the marginal statement meaning that there exists a
residually solvable linear group  which does not embed in its Malcev
completion. By  residually solvable we mean a group that embeds in its
completion with respect to  all finitely generated solvable representations.

The following example came out from a discussion with D. Arapura, J\'anos
Koll\'ar, M. Nori, T. Pantev,  M. Ramachandran and D. Toledo.

Consider nontrivial  smooth family of  smooth abelian varieties of dimension
$N$ over curve $C$. Let us denote this family by $X$. The fundamental group
of $X$ is given by the following exact sequence
\[ 0 \longrightarrow {\Bbb Z}^{2N} \longrightarrow \pi_{1}(X) \longrightarrow
\pi_{1}(C) \longrightarrow 0.\]
The group $\pi_{1}(X) $ is a semidirect product of  the groups ${\Bbb Z}^{2N}$
and $\pi_{1}(C)$. For generic enough family we can make the monodromy action
\[M:  \pi_{1}(C) \longrightarrow Sp(2N, {\Bbb Z}) \]
to be irreducible and from here one can get that the image of ${\Bbb Z}^{2N}$
in $H^{1}(X,{\Bbb Z})$  is trivial.

The the group  $\pi_{1}(X)$ is linear. To see that we consider the morphism
\[l:  \pi_{1}(X) \longrightarrow  SL(2,{\Bbb C}) \times [ GL(V) \ltimes V ].\]
Here $V$ is a vector space over $\Bbb{C}$ of dimension $N$ on which
${\Bbb Z}^{2N}$ acts discretely. It is easy to see that $L$ is an  injection.

The group $\pi_{1}(X)$ is also virtually residually solvable. This can be
seen as follows:

Choose a prime number $p$. Since the group $\pi_{1}(X)$ is linear, namely it
embeds in $GL(T)$  for some vector space $T$ we can embed it up to a finite
index  in a series of finite solvable groups  $GL(T_{p^{q}})$  for $q=1,2,
\ldots$.

{}From another point $\pi_{1}(X)$ does not embed in its Malcev completion up
to a finite index.
Assume that  $\pi_{1}(X)$ does  embed in its Malcev completion. Then
lemma 2.3 and lemma 2.4  imply that if the image of ${\Bbb Z}^{2N}$ in
$H^{1}(X,{\Bbb Z})$  is trivial then the image of ${\Bbb Z}^{2N}$ in the Malcev
completion of $\pi_{1}(X)$ is trivial which is not the case in our situation.
Therefore our technique does not answer the question if the universal
coverings of $X$ or of generic hyperplane sections of it is holomorphically
convex. Of course this is true and can be seen as follows:

\begin{prop} The universal covering of  any  smooth family of Abelian varieties
or of any    generic hyperplane section of  them  is holomorphically  convex.

\end{prop}

{\bf Proof.} It follows from  \cite{K1}, Theorem 6.3 that every smooth
family of abelian varieties over a curve has a linear fundamental group since
according
to 6.3 \cite{K1} after a  finite etale  covering it is birational to a   family
of a smooth abelian varieties. But the universal covering of a  family of
smooth abelian varieties or a generic hyperplane section of it is
holomorphically convex since it is a Stein space since. It embeds in $(SIEG
\times C^{N})$, where $SIEG $ is the Siegel upperhalf plane.

\hfill $\Box$

 It also follows from  \cite{LM} where  more
powerful technique, the theory of  Nonabelian Mixed Hodge Structures, is used.

We  formulate:

\begin{corr} Let $X$ be a smooth projective variety    with an infinite
virtually nilpotent   fundamental group and such  that ${\rm rank}\,
{\rm Pic}(X) = 1$ (or better ${\rm rank}\, NS(X)=1$). Then for every
subvariety $Z$
in $X$ we have that  $ \im [\pi_{1}(Z)\longrightarrow \pi_{1}(X)]$ is infinite.
 \end{corr}

The proof is an easy consequence of  \cite{K2} (Chapter 1).

We demonstrate a quick application of the above corollary. Denote by $A$ a
four
dimensional abelian variety. Let us say that $X$ is a hypersurface in it with
an isolated singular point $s$ of the following type $xy=zt$. Denote by $X'$
the blow of $X$ in $s$.
We glue in $X$ ${\Bbb P}^{1} \times {\Bbb P}^{1}$ instead of $s$.
It is easy to see that $X'$ is smooth and that ${\Bbb P}^{1}$ can be blown
down. The new space $X''$, obtained after blowing down one of the above
${\Bbb P}^{1}$, is smooth, ${\rm rank} NS(X'')=1$ and
$\pi_{1}(X)={\Bbb Z}^{4}$. Therefore the Shafarevich conjecture should be true
for $X''$. But as we can see $X''$ contains ${\Bbb P}^{1}$. This contradicts
the
above corollary and we conclude that $X''$ is not projective. Of course all
this can be seen in many different ways. This is another spelling of the fact
that ${\Bbb P}^{1}$, that remains in $X''$, should be homologically nontrivial
if $X''$ is projective.

We give now an idea of an alternative proof of Theorem 1.1 which came from
conversations with M. Ramachandran. It is based on the use of $\pi_{1}(X)$
equivariant harmonic maps to the universal coverings to Higher Albanese
varieties defined in \cite{HZ}. Combined with the strictness property for the
nonabelian Hodge theory this seems to be a very promising idea ( see
\cite{LM}).

Denote by $G_{s}$ the complex simply connected group $\pi_{1}(X) /
\Gamma^{s+1}$, where $\Gamma^{i}$ are the groups from the lower central series
for $\pi_{1}(X)$ and $\Gamma^{s}$ is the smallest nontrivial one. The
corresponding Lie algebra $g_{s}$ has MHS. Denote by $F^{0}G_{s}$ the closed
subgroup in $G_{s}$ group that corresponds to  $F^{0}g_{s}$. Since the group
$\pi_{1}(X) / \Gamma^{s+1}$ is unipotent then as it is easy to see we have a
free action of  the corresponding to $G_{s}$ lattice $G_{s}({\Bbb Z})$ on
$G_{s}/F^{0}G_{s}$.

Therefore in the same way as in \cite{KR} we obtain a $\pi_{1}(X) $
equivariant proper horizontal holomorphic map ( see \cite{HZ})
\[ \widetilde{X} \longrightarrow G_{s}/F^{0}G_{s}.\]
According to \cite{H1} $G_{s}/F^{0}G_{s}$ is biholomorphic to $\Bbb{C}^{N}$.
Therefore $\widetilde{X}$ is holomorphically convex.

\begin{rem}{\rm The above argument is weaker then the argument we have used in
the first proof. It cannot be generalized to the case of residually nilpotent
groups since in this case $G_{s}/F^{0}G_{s}$ will not be a manifold.}
\end{rem}

\section{Solvable coverings}

We would like to obtain the analog of Theorem 1.1 for solvable
coverings. The analog of Theorem 1.2 for solvable groups  - Corollary 2.3
was proved in the  previous section as a consequence of the result of Arapura
and Nori. We cannot prove solvable analog of theorem 1.1. The maximum we can do
is to realize how close the solvable representations come to nilpotent ones. To
be able   to do so we need  to generalize slightly the result of  Arapura and
Nori.

First we  prove Theorem 1.3.

\noindent
{\bf Proof.} (The idea of the proof was suggested to me by  T. Pantev.) Denote
by $\Gamma$  the image of the solvable representation $\rho : \pi_{1}(X)
\rightarrow L$.  We need to show that  either $\Gamma$ is virtually nilpotent
or there exists a holomorphic map with connected fibers $f : X \rightarrow C$
to a smooth curve $C$ of genus $\geq 1$.

First  we  introduce some notations. For a
finitely generated group $\Gamma$ denote by $\Sigma(\Gamma)$ the set of all
special characters of $\Gamma$. That is
\[
\Sigma(\Gamma) := \left\{ \alpha : \Gamma \rightarrow {\Bbb C}^{\times}
\left| H^{1}(\Gamma, {\Bbb C}_{\alpha}) \neq 0 \right. \right\},
\]
where ${\Bbb C}_{\alpha}$ is the one dimensional $\Gamma$-module associated to
$\alpha$. Now  we recall the following:

\begin{prop} [Arapura-Nori \cite{AN}]  Let $\Gamma$ be a finitely generated
${\Bbb Q}$-linear solvable group. Then the following are equivalent
\begin{enumerate}
\item $\Gamma$ is virtually nilpotent.
\item $\Sigma(\Gamma)$ consists of finitely many torsion characters.
\end{enumerate}
\end{prop}

Due to this proposition it is enough to show that either $\Sigma(\Gamma)$
consists of finitely many torsion characters or $X$ has a non-trivial map to
a curve of genus bigger than zero.  Now, since $\pi_{1}(X)$ surjects on
$\Gamma$ it follows that $\Sigma(\Gamma) \subset \Sigma(\pi_{1}(X))$ and hence
it suffices to show that either $\Sigma(\pi_{1}(X))$ consists of finitely
many torsion characters or $X$ has an irrational pencil.

\medskip

For a smooth projective variety $X$ denote by $M(X)$ the moduli space of
homomorphisms from $\pi_{1}(X)$ to ${\Bbb C}^{\times}$. The locus of special
characters is a jump locus in $M(X)$ and hence it is a  subscheme in a natural
way.
It turns out that $\Sigma(\pi_{1}(X))$ is actually a smoooth subvariety having
very special geometric properties which we are going to exploit. Since the
subvariety $\Sigma(\pi_{1}(X)) \subset M(X)$ is completely canonical one
expects it to have an intrinsic description. One way to construct natural
subvarieties in $M(X)$ is via pullbacks. Namely, given any surjective morphism
$\varphi : X \rightarrow Y$ we can pullback the moduli space of characters of
$\pi_{1}(Y)$ to get a subvariety $\varphi^{*}M(Y) \subset M(X)$. According to
\cite{SA}, Lemma 2.1 and Theorem 6.1 every connected component $\Sigma$ of the
subvariety $\Sigma(\pi_{1}(X)) \subset M(X)$ is of this kind. More specifically
for every such $\Sigma$ there exists a torsion character $\sigma \in \Sigma$
and a connected abelian subvariety $P \subset \Alb(X)$ so that $\Sigma$ is
the translation of $\varphi^{*}M(\Alb(X)/P) \subset M(X)$ by $\sigma$. Here
$\varphi : X
\rightarrow \Alb(X) \rightarrow \Alb(X)/P$ is the composition of the Albanese
map and the natural quotient morphism. In particular, $\Sigma(\pi_{1}(X))$
has a positive dimensional component if and only if its intersection with the
set of all unitary characters has a positive dimensional component. Now the
Hodge decomposition of the cohomology of a unitary local system implies that
unless $\Sigma(\pi_{1}(X))$ consists of finitely many torsion characters
the subvariety of all special line bundles in ${\rm Pic}^{\tau}(X)$ has a
positive dimensional component. Indeed, for a unitary character $\alpha$
denote by ${\Bbb L}_{\alpha}$ the corresponding rank one local system and
by $L_{\alpha} = {\Bbb L}_{\alpha}\otimes_{\Bbb C} {\cal O}_{X}$ the
corresponding holomorphic line bundle. Now by the Hodge theorem
\[h^{1}(\pi_{1}(X),{\Bbb C}_{\alpha})
= h^{1}(X, {\Bbb L}_{\alpha}) = h^{1}(X,L_{\alpha})+h^{0}(X,\Omega^{1}_{X}
\otimes L_{\alpha}) = 2h^{1}(X,L_{\alpha}),\]
i.e. $\alpha$ is a special character iff the line bundle $L_{\alpha}$ is
special.

Furthermore a theorem of Beauville  (\cite{BE}, Proposition 1) asserts that
the subvariety of ${\rm Pic}^{0}(X)$ consisting of special line bundles is a
union of a finite set and the subvarieties of the form $f^{*}{\rm Pic}^{0}(B)$
where $f : X \rightarrow B$ is a morphism with connected fibers to a curve $B$
of genus $\geq 1$. Thus $X$ posseses irrational pencils which finishes the
proof of Theorem 1.3 \hfill $\Box$

\medskip

The above theorem can be seen as the solvable analog of the theorem of
Simpson's that $SL(n,{\Bbb Z})$ is not a K\"{a}hler group, $n>2$. This
theorem gives a  way of constructing new examples of non-K\"{a}hler groups.
In particular any group $\Gamma$ with infinite $H^{1}([\Gamma,\Gamma],
{\Bbb Q})$ possessing a solvable linear quotient defined over ${\Bbb
Q}$ that is not virtually nilpotent cannot be K\"{a}hler.

\medskip

Now we prove theorem 1.4.

\noindent
{\bf Proof.} We would like to use theorem 3.1. Therefore we need an infinite
solvable representation defined over ${\Bbb
Q}$. First we show

\begin{lemma} Let $S$ be a connected affine solvable group defined over ${\Bbb
Q}$. If $\rho : \pi_{1}(X) \rightarrow S({\Bbb C})$ is a Zariski dense
represenattion into the group of complex points of $S$, then $\rho$ can be
deformed to a representation $\nu : \pi_{1}(X) \rightarrow S(\overline{{\Bbb
Q}})$ having an infinite image.
\end{lemma}
{\bf Proof.}
Denote by $\Lambda$ the image of $\rho:\pi_{1}(X) \to S$. Let $\phi:\Lambda
\to S(\Bbb{C})$ be a homomorphism from $\Lambda$
to the group of complex points of $S$ with an  infinite image.
We want to find a subgroup $\Lambda'$ of $ \Lambda$ of finite index and a
homomorphism $\phi':\Lambda'\to S (\overline{{\Bbb Q}})$ with infinite image,
such that $\phi '$ is arbitrarily close to $\phi$.

Let $S_{0} := S$, and consider the (upper) derived series $S_{i+1} :=
[S_{i},S_{i}]$ for $S$.
We choose the maximal $i$ such that $\Lambda\cap S_{i(\Bbb{C})}$ is of finite
index of $\Lambda$.
We replace $G$ by $\Lambda' = \Lambda\cap S_i(\Bbb{C})$.  The image of
$\Lambda'$ in $S_i(\Bbb{C})/S_{i+1}(\Bbb{C})$
is infinite.  The group $A = S_{i}/S_{i+1}$ is either a torus or a vector
space group.  Let $X$ be the affine variety of homomorphisms $\Lambda' \to A$,
$Y$ the affine variety  $Hom(\Lambda',S_{i}), X' $ the image of $Y$ in $X$.
Thus $X'$ is the affine subvariety of $X$ consisting of homomorphisms which
factor through $S_{i}$.  We want to find points on $X'(\overline{{\Bbb Q}})$
arbitarily close in $X'(\Bbb{C})$  to the point  defined by $\phi$.  If the
original point is defined
over $\overline{{\Bbb Q}}$, we are done.  If not, the original point cannot
be isolated,
since $X'$ is defined over $\overline{{\Bbb Q}}$.  Thus we have arbitrarily
close points
defined over $\overline{{\Bbb Q}}$.  To find a $\overline{{\Bbb Q}}$ point
with an infinite image we consider two cases:

1) $A$ is a vector space group. Let $g$ be an element of
$\Lambda'$ such that $\phi(g) $ maps to a non-trivial element of $A$. Then for
every nearby representation $\phi'$ the element $\phi'(g)$ is non-trivial
in $A$, therefore it is of infinite order in $A$.

2) $A$ is a torus $T$.
We fix an element $g$ in $\Lambda'$ such that $\phi(g)$ maps to a point of
infinite order on $T$.  Let $Z \subset T$ denote the image of $X'$ in $T$
under the map which takes each homomorphism $\Lambda'\to T$ to the
image of $g$ in $T$.  By definition $X' \to Z$ is a surjective map and it is
defined over
$\overline{{\Bbb Q}}$, so every $\overline{{\Bbb Q}}$ -point of $Z$ comes from
a $\overline{{\Bbb Q}}$-point of $X'$
and therefore from  a $\overline{{\Bbb Q}}$ -point of $Y$, i.e. an actual $\bar
Q$ -homomorphism
from $\Lambda'$  to $S_{i}$.  So it is  enough to find a $\overline{{\Bbb Q}}$
-point of  $Z$  which is
close to the image of the original homomorphism $\phi$  but which is also
of infinite order.

We prove the following:

\begin{claim} Let $T$ be a torus, $Z$ a $\overline{{\Bbb Q}}$-affine subvariety
of $T$, $p$ a point in
$Z(C) $ of infinite order.  Then $p$ is in the closure of the subset of
$Z(\overline{{\Bbb Q}}) $ consisting of points of infinite order.
\end{claim}
{\bf Proof.} If $p$  is in $Z(\overline{{\Bbb Q}})$, we  are done.    If not,
there exists a character
$\chi: Z \to GL(1,\Bbb{C})$ such that $\chi(p)$ is not in $ \overline{{\Bbb
Q}}$ .  As $Z, T$  and $\chi$
are defined over $ \overline{{\Bbb Q}}$ , and $GL(1,\Bbb{C})$  is
1-dimensional, it follows
that $\chi(Z)$  is an open subset of $GL(1,\Bbb{C})$.  In particular, every
neighborhood of $p$ in $Z$ maps to a neighborhood of $\chi(p)$ containing
non-torsion elements.  If $ q \in Z(\overline{{\Bbb Q}})$ maps to a non-torsion
element,
then of course $q$ is a point of infinite order in $T(\Bbb{C})$. This finishes
the proof of the claim and the lemma.
\hfill $\Box$

\medskip
To get  an infinite solvable representation defined over ${\Bbb Q}$
consider the affine solvable group $\widetilde{S}$ obtained from $S$ by
restriction of scalars, i.e. $\widetilde{S} := {\rm res}_{\overline{{\Bbb Q}}/
{\Bbb Q}}S$. The representation $\nu$ induces a representation $\tilde{\nu}:
\pi_{1}(X) \rightarrow \widetilde{S}({\Bbb Q})$ which has an image isomorphic
to the image of $\nu$. \hfill $\Box$

\hfill $\Box$

Now to prove  Corollary~1.1 we have to consider the following two alternatives

1) The group $\tilde{\nu}(\pi_{1}(X)$
is virtually nilpotent  so it has a subgroup of
finite index which is nilpotent.

2) There exists an holomorphic map with connected fibers $f : X
\rightarrow C$ to a curve of genus one or higher.  But then we know that $
\pi_{1}(C)$ embeds in its Malcev completion.

 In both cases  there exists a finite \'{e}tale cover of $X$ which has a
non-trivial Albanese variety, which is what we need.

Theorem 1.5 follows easily from Theorem 1.1 and the result from \cite{KR}.

{\bf Proof.} Let us start with a complex linear representation
${\rm im}[\pi_{1}(X) \rightarrow L]$. Then we have the following three
possibilities.

(a) The image of $\pi_{1}(X)$ in $L/R^{u}L$ does not have zero or two ends.
Then we can apply \cite{KR}   to conclude that $\widetilde{X}$ has
a non-constant holomorphic function.

(b) The image ${\rm im}[\pi_{1}(X) \rightarrow L/R^{u}L]$ has two ends. Then
by the theorem of Hopf and Freudenthal it follows that ${\rm im}[\pi_{1}(X)
\rightarrow L/R^{u}L]$ has a subgroup of finite index that is isomorphic to
${\Bbb Z}$. Therefore the abealianization of $\pi_{1}(X)$ is not finite. This
implies that the Malcev representation is not trivial and we apply
Theorem 1.1 to finish the proof.

(c) The image of $\pi_{1}(X)$ in $L/R^{u}L$ has zero ends. Then $L/R^{u}L$ is
a finite group. So the Malcev representation is not trivial and we are taken
applying Theorem 1.1.
\hfill $\Box$

\medskip

Theorem 1.5  follows from   \cite{LM} as well.  To be able to attack conjecture
1.1 we should be able to analyze the real issue,  the semisimple
representations. Some initial steps in this  direction are done in \cite{LM}.

\medskip

What should we do  if the answer of {\bf Question 3} is positive? We hope
using \cite{LP}  to be able to handle the case when the image of $ \pi_{1}(X)$
  in its proalgebraic completion is infinite.

\medskip

What should we do  if the answer of {\bf Question 2} is positive and  the
$ \pi_{1}(X)$  in  question has finite image in its proalgebraic completion.
At the moment this case seems to be out of reach.

\medskip

\noindent
MSRI, 1000 Cent. Drive, Berkeley, CA 94720


\begin{thebibliography}{40}
\bibitem[1]{A}{\bf  D. Arapura, } {\em Higgs line bundles Green Lazarsfeld
sets and maps   K\"{a}hler manifolds to curves.\/} Bull. of A.M.S.,Vol. 26,
1992, pp. 310-314.
\bibitem[2]{AN}{\bf  D. Arapura, M. Nori} {\em  Solvable  fundamental groups
of algebraic vatieties and  K\"{a}hler manifolds.\/} In preparation.
\bibitem[3]{AR}{\bf  D. Arapura, P. Bressler, M. Ramachandran} {\em  On the
fundamental group of a compact K\"{a}hler manifold .\/} Duke Mathematical
Journal, v.68, No.3, (1992), 477-488.
\bibitem[4]{BU}{\bf  G. Baumshlag} {\em Topics in combinatorial group
theory .\/} Lectures in mathematics, ETH, Zurich, Birkhauser, 1993.
\bibitem[5]{BE}{\bf  A. Beauville} {\em Annulation de $H^{1}$ et systemes
paracanonidues sur les surfaces .\/} Crelle Journ., Vol. 388, 1988, pp.
149-157.
\bibitem[6]{BL}{\bf   F. Bogomolov, L. Katzarkov} {\em  Projective surfaces
with intersting fundamental groups.\/} in preparation , (1995).
\bibitem[7]{CM}{\bf  F. Campana }{\em  Remarques sur le revetement universal
des varietes  K\"{a}hleriennes  compactes.\/} Bull. Soc. Math. France, 122,
1994, pp. 255-284.
\bibitem[8]{CM1}{\bf  F. Campana }{\em  Remarques sur le groupes de
K\"{a}hleriennes  nilpotent.\/} Ann. Scien. L'Ecole Norm. Sup. 28, fasc.3,
1995, pp. 307-316.
\bibitem[9]{CT}{\bf  J. Carlson, D. Toledo} {\em Quadratic Presentations and
Nilpotent K\"{a}hler Groups.\/} Preprint.
\bibitem[10]{C}{\bf  K. Corlette} {\em Flat $G$-bundles with canonical
metrics.\/} J. Differential Geometry, v. 28, (1988), pp. 361-382.
\bibitem[11]{D}{\bf P. Deligne } {\em Theorie de Hodge II.\/} Publ. Math.
I.H.E.S., 40, (1971), pp. 5-58.
\bibitem[12]{DG}{\bf  P. Deligne, P. Griffiths, J. Morgan, D. Sullivan} {\em
Real Homotopy theory of  K\"{a}hler manifolds  .\/} Inventiones Matematicae,
29, (1975), pp. 245-274.
\bibitem[13]{GM}{\bf M. Goresky, R. MacPherson} {\em Stratified Morse Theory.
\/} Springer Verlag, 1988.
\bibitem[14]{GL}{\bf M. Green, R. Lazarsfeld} {\em Higher obstructions to
deforming cohomology groups of line bundles  \/} J. Amer. Math. Soc., {\bf  4},
87-103.
\bibitem[15]{GS}{\bf  M. Gromov , R. Schoen} {\em Harmonic maps into singular
spaces and p-adic superrigidity  for lattices in groups of rank one .\/} IHES -
 Publications Mathematiques, 76, 1992, pp. 165-246.
\bibitem[16]{H}{\bf  R. Hain}  {\em The de Rham homotopy theory of complex
algebraic variety I .\/}  K-Theory, 1, 1987, pp. 271-324.
\bibitem[17]{H1}{\bf R. Hain and S. Zucker,} {\it Unipotent variations of
mixed Hodge structures,\/} Invent. math. {\bf 88} (1987) 83-124.
\bibitem[18]{HZ}{\bf R. Hain and S. Zucker,} {\it A guide to Unipotent
variations of
mixed Hodge structures,\/} Hodge theory, Springer LNM, 1246, 1987,pp. 92-106.
\bibitem[19]{HI}{\bf N. J. Hitchin} {\em The self-duality equations on a
Riemann surface.\/} Proc. Lond.  Math. Soc., (3) vol.55, 1987, 59-126.
\bibitem[20]{BR}{\bf  B. Lasell,  M. Ramachandran} {\em Observations on
harmonic maps and singular vbarieties.\/}  To appear in Ann. Ec. Norm. Sup..
\bibitem[21]{LM}{\bf L.Katzarkov } {\em On the Shafarevich maps.\/} Preprint,
1995.
\bibitem[22]{KR}{\bf L.Katzarkov M. Ramachandran} {\em On the Shafarevich
conjecture for surfaces.\/} Preprint, 1995.
\bibitem[23]{LP}{\bf L.Katzarkov T. Pantev } {\em Nonabelian Mixed Hodge
Fiction.\/} Preprint, 1995.
\bibitem[24]{K1}{\bf  J\'anos Koll\'ar} {\em Shafarevich maps and plurigenera
of algebraic varieties  .\/} Inventiones Matematicae, 113, Fasc. 1, 1993, pp.
165-215.
\bibitem[25]{K2}{\bf  J\'anos Koll\'ar} {\em Shafarevich maps and automorphic
forms  .\/} To appear in Princeton University Lecture Notes.
\bibitem[26]{N1}{\bf  T. Napier }{\em  Convexity properties of coverings of
smooth projective varieties.\/}  Mathematische Annalen, 286, 1990, pp. 433-480.
\bibitem[27]{SA}{\bf C. Simpson} {\em Subspaces of moduli spaces of rank one
local systems. \/} Ann. Scien. Ecol. Norm. Sup., t. 26 , 1993, pp. 361-401.
\bibitem[28]{SC}{\bf  C. T. Simpson} {\em Constructing variations of Hodge
structures  using Yang-Mills theory and applications to uniformization.\/}
Journal of A.M.S, 1 (1988),  pp. 368- 391.
\bibitem[29]{SL}{\bf  C. T. Simpson} {\em Higgs bundles and local systems.\/}
Publ.  Math. I.H.E.S, 75, 1992, pp. 5-95.
\bibitem[30]{SIM1}{\bf  C. T. Simpson} {\em Moduli spaces of representations 1.
\/} IHES -  Publications Mathematiques, 79, (1995), pp. 5-79.
\bibitem[31]{SIM2}{\bf  C. T. Simpson} {\em Moduli spaces of representations 2.
\/} IHES -  Publications Mathematiques, 80, (1995), pp. 5-79.
\bibitem[32]{SIM3}{\bf  C. T. Simpson} {\em Nonabelian Hodge theory \/}
International congress of Mathematics  Kyoto, 1990, Proceedings, Springer -
Tokyo, (1991), pp. 747-756.
\bibitem[33]{SV}{\bf  A. J. Sommese, A. Van de Ven }{\em  Homotopy groups of
pullbacks of varieties.\/} Nagoja Math. Journal , vol. 102, 1986, pp. 79-90.


\end{thebibliography}
\end{document}